\documentclass[10pt,reqno]{amsart}
\usepackage{geometry}
\usepackage{amsmath}
\begin{document}
\author{Tatenda Kubalalika}     
\email{tatendakubalalika@yahoo.com}
\title{The Riemann Hypothesis is false}
\maketitle
\textbf{ABSTRACT.} Let $\Theta$ denote the supremum of the real parts of the zeros of the Riemann zeta function $\zeta(s)$. By carefully exploiting the Landau-Gonek explicit formula, the Dirichlet polynomial approximation of $\zeta(s)$ and other classical results in analytic number theory, we rigorously demonstrate that $\Theta=1$. This entails the existence of infinitely many Riemann zeros off the critical line (thus disproving the Riemann Hypothesis (RH)). The paper is concluded by a brief discussion of why our argument does not work for both Weil and Beurling zeta functions, whose analogues of the RH are known to be true. 
\\
\\
\textit{2020 Mathematics Subject Classifications: 11M26, 11M06.}
\\
\\
\footnote{\textbf{Data availability statement:} This manuscript has no associated data.}
\textbf{Introduction.} 
The Riemann zeta function is a function of the complex variable $s$, defined in the half-plane $\Re(s)>1$ by $\zeta(s):=\sum_{n=1}^\infty n^{-s}$ and in the whole complex plane by holomorphic continuation. Euler noticed that for $\Re(s)>1$, $\zeta(s)$ can be expressed as a product $\prod_{p}(1-p^{-s})^{-1}$ over the entire set of primes, which entails that $\zeta(s)\neq 0$ for $\Re(s)>1$. It can be shown that $\zeta(s)$ extends to $\mathbb{C}$ as a meromorphic function with only a simple pole at $s=1$, with residue $1$, and satisfies the functional equation $\xi(s)=\xi(1-s)$, where $\xi(z) = \frac{1}{2}z(z-1)\pi^{-z/2}\Gamma(\frac{1}{2}z)\zeta(z)$ and $\Gamma(w) = \int_{0}^{\infty} e^{-x}x^{w-1} \mathrm{d}x$. Let $\chi$ be a primitive Dirichlet character and $L(s, \chi)$ be the corresponding L-function. Most of the analytic properties for $\zeta$ are also known to hold for $L(s, \chi)$ [5, Chapter 10]. Let $\Theta_{\chi}$ be the supremum of the real parts of the zeros of $L(s, \chi)$. The Generalised Riemann Hypothesis (GRH) is equivalent to the statement that $\Theta_{\chi}=\frac{1}{2}$ for each $\chi$, and the Generalised Prime Number Theorem is equivalent to the statement that $L(s, \chi) \neq 0$ at $\Re(s)=1$ for each $\chi$. For a far more thorough discussion of the RH and GRH, the interested reader is kindly referred to Borwein \textit{et al} [1]. 
\\
\\

\textbf{Definitions.} Let: $\rho$ denote a non-trivial zero of $\zeta, \sigma=\Re(s)$, $\gamma=\Im(\rho),\Lambda$ be the von-Mangoldt function and $\Theta \in [\frac{1}{2}, 1]$ be the supremum of the real parts of the $\rho$'s. 
\\
\\
\textbf{Lemma 1. (Theorem 4.11 of [7]).} \textit{We have $\zeta(s)=\sum_{n \leq x} n^{-s} -\frac{x^{1-s}}{1-s} + O(x^{-\sigma})$ uniformly for $\sigma\geq \sigma_{0}>0$ and $|\Im(s)| < 2\pi x/C$, where $C$ is a given constant $>1$.}
\\
\\
\textbf{Lemma 2, (Theorem 1 of [3]).} \textit{One has $\sum_{0 \leq \gamma \leq  y}  n^{\rho-1}  =-\frac{y\Lambda(n)}{2\pi n}  + O((\log y)\log \log 3n)$ uniformly for $2 \leq n \ll y.$}
\\
\\
\textbf{Theorem 1.} \textit{One has $\Theta=1$.} 
\\
\begin{proof} 
Let: $T \in \mathbb{R}_{>0}$ be arbitrarily large and $X \geq 2$. By Lemma 12.2 of [5], we know that there exists some $X_0 \in [X, X+1]$ such that $\frac{\zeta'}{\zeta}(\sigma + iX_0) \ll \log^{2} X$ uniformly for $-1\leq \sigma \leq 2$. Thus we can choose $T$ and $H \in [2T, 2T+1]$ such that 
\begin{equation}
    \max\Big(\Big|\frac{\zeta'}{\zeta}(\sigma+iT)\Big|, \Big|\frac{\zeta'}{\zeta}(\sigma+iH)\Big| \Big) \ll \log^{2}T
\end{equation}
uniformly for 
$-1 \leq \sigma \leq 2$. By the residue theorem, we have 
\begin{align}
    \sum_{T \leq \gamma_{T} \leq H}(2n)^{\rho_{T}-1} &= \frac{1}{2\pi i}\Bigg(\int_{1+\frac{1}{\log 2n}+iT}^{1+\frac{1}{\log 2n}+iH} + \int_{1+\frac{1}{\log 2n}+iH}^{-1+iH} +\int_{-1+iH}^{-1+iT} + \int_{-1+iT}^{1+\frac{1}{\log 2n}+iT}  \Bigg) \frac{\zeta'}{\zeta}(s)(2n)^{s-1} \mathrm{d}s \\
    &=I_1 + I_2 + I_3 + I_4,
\end{align}
say. From now on, let $T \leq t \leq H.$
By the approximate functional equation for $\zeta$ (Theorem 4.15 of [7]), note that 
\begin{equation}
    \zeta(\sigma-it)=\sum_{n \leq x} n^{-\sigma+it} +\chi(\sigma-it)\sum_{v \leq y} v^{\sigma-1-it} + O(x^{-\sigma} + t^{\frac{1}{2}-\sigma}y^{\sigma-1})
\end{equation}
where $0\leq \sigma \leq 1, x>h>0, y>h>0, t=2\pi xy, h>0$ is a constant, $\chi(s)=\frac{\zeta(s)}{\zeta(1-s)}$ and $|\chi(\sigma-it)|=|\chi(\sigma+it)|=\Big| \Big(\frac{t}{2\pi}\Big)^{\frac{1}{2}-\sigma-it}e^{i(t+\frac{\pi}{4})}(1+O(1/t))\Big|\sim(t/2\pi)^{\frac{1}{2}-\sigma}$ [7, p.78]. Thus since $|\zeta(-it)|\asymp t^{1/2}|\zeta(1+it)|\ll T^{1/2}\log T$, it follows by taking $\sigma=0$ in (4) that 
\begin{equation}
    S(x):=\sum_{n \leq x} (2n)^{it} \ll T^{\frac{1}{2}}\log^{3} T 
\end{equation} 
uniformly for $\frac{T}{2\log^{2}T} \leq x \leq T$. Let $f(n)=(2n)^{\sigma-1}.$ Then by Abel summation, it follows that 
\begin{equation}
    \sum_{\frac{T}{\log^{2} T} \leq 2n \leq T} (2n)^{\sigma-1+it} = S(T/2)f(T/2)-S(T/2\log^{2}T)f(T/2\log^{2}T) - \int_{T/2\log^{2}T}^{T/2} S(x)f'(x) \mathrm{d}x \ll \frac{\log^{7} T}{T^{\frac{1}{2}-\sigma}}
\end{equation}
for $-1 \leq \sigma \leq 2.$ Combining (1) and (6) gives $\sum_{\frac{T}{\log^{2}T} \leq 2n \leq T} I_{2} + I_{4} \ll T^{1/2}\log^{8}T.$ By the log-derivative of the functional equation for $\zeta$, we have 
\begin{equation}
    \frac{\zeta'}{\zeta}(-\sigma+it) \ll \log T
\end{equation} 
for $\sigma \geq \sigma_0 >0$. Combining (6) and (7) yields $\sum_{\frac{T}{\log^{2}T} \leq 2n \leq T} I_{3} \ll T^{1/2}\log^{8}T.$ Thus we arrive at 
\begin{equation}
    \sum_{\frac{T}{\log^{2}T} \leq 2n \leq T} I_{2} + I_{3} + I_{4} \ll T^{1/2}\log^{8} T.
\end{equation}
Note that $\frac{\zeta'}{\zeta}(s)=-\sum_{m=2}^{\infty} \frac{\Lambda(m)}{m^s}$ for $\sigma>1$. Inserting this into $I_1$ gives 
\begin{align}
I_1 &= \frac{(T-H)\Lambda(2n)}{4\pi n}  + \frac{e}{2\pi i}\sum_{m=2,m \neq 2n}^{\infty} \Bigg( \frac{\Lambda(m)(2n)^{iT}}{m^{1+\frac{1}{\log 2n}+iT}\log(2n/m)} - \frac{\Lambda(m)(2n)^{iH}}{m^{1+\frac{1}{\log 2n}+iH}\log(2n/m)} \Bigg). 
\end{align}
 From now on, let $\varepsilon=\varepsilon_{T}=T^{-1/4}$ and $\frac{T}{\log^{2} T} \leq 2n \leq T$ (unless specified otherwise). For $2 \leq m \leq \frac{T}{(1+\varepsilon)\log^{2} T}$  and $m \geq T(1+\varepsilon)$, note that $\frac{1}{|\log(2n/m)|}$ is uniformly $\ll 1/\varepsilon=T^{1/4}$. For $m \neq 2n$, let $g(n)=\frac{1}{m^{1+\frac{1}{\log 2n}}\log (2n/m)}$ hence $g'(n)=\frac{m^{-1-\frac{1}{\log 2n}}(\log(2n)\log(m)-\log^{2}(2n)-\log^{2}(m))}{n\log^{2}(2n)\log^{2}(2n/m)}.$ Note that $g'(n)$ is strictly negative, hence $|g'(n)|=-g'(n).$ Thus in similar fashion to (6), we obtain from (5) that 
\begin{align}
    \sum_{\frac{T}{\log^{2} T}\leq 2n \leq T} \frac{(2n)^{it}}{m^{1+\frac{1}{\log 2n}}\log(2n/m)} &= S(T/2)g(T/2) - S(T/2\log^{2} T)g(T/2\log^{2} T) - \int_{T/2\log^{2} T}^{T/2} S(x)g'(x) \mathrm{d}x \\
    &\ll \frac{T^{3/4}\log^{5} T}{m^{1+\frac{1}{\log T}}} - T^{1/2}\log^{3}T\int_{T/2\log^{2} T}^{T/2} g'(x) \mathrm{d}x \\
    & \ll \frac{T^{3/4}\log^{5} T}{m^{1+\frac{1}{\log T}}}
\end{align}
uniformly for $2\leq m \leq \frac{T}{(1+\varepsilon)\log^{2}T}$ and $m \geq T(1+\varepsilon).$ For $2n \leq \frac{m}{1+\varepsilon}$ and $2n \geq m(1+\varepsilon)$, $\frac{1}{|\log(2n/m)|}$ is also uniformly $\ll 1/\varepsilon = T^{1/4}$ so it similarly follows from (5) by Abel summation that 
\begin{equation}
    \sum_{\frac{T}{\log^{2} T}\leq 2n \leq \frac{m}{1+\varepsilon} \leq T} \frac{(2n)^{it}}{m^{1+\frac{1}{\log 2n}}\log(2n/m)} \ll \frac{T^{3/4}\log^{5} T}{m^{1+\frac{1}{\log T}}} 
\end{equation} and 
\begin{equation}
    \sum_{\frac{T}{\log^{2} T} \leq m(1+\varepsilon)\leq 2n \leq T} \frac{(2n)^{it}}{m^{1+\frac{1}{\log 2n}}\log(2n/m)} \ll \frac{T^{3/4}\log^{5} T}{m^{1+\frac{1}{\log T}}}.
\end{equation}
We now treat the ranges $m+1 \leq 2n \leq m(1+\varepsilon)$ and $\frac{m}{1+\varepsilon}\leq 2n \leq m-1$, where $\frac{T}{(1+\varepsilon)\log^{2}T} \leq m \leq T(1+\varepsilon)$. For $m+1 \leq 2n \leq m(1+\varepsilon)$, put $2n=m+r$ where $1\leq r \leq m\varepsilon-1$ thus $\log(2n/m)=\frac{r}{m}(1 + O(r/m))$ hence $\frac{1}{\log(2n/m)}=\frac{m}{r}(1+O(r/m))=\frac{m}{r}+O(1)$ hence
\begin{align}
    \sum_{m+1 \leq 2n \leq m(1+\varepsilon)} \frac{\Lambda(m)(2n)^{iT}}{m^{1+\frac{1}{\log 2n}+iT}\log(2n/m)} & = \frac{\Lambda(m)}{m^{iT}}\sum_{m+1 \leq 2n \leq m(1+\varepsilon)} \frac{(2n)^{iT}}{m^{1+\frac{1}{\log 2n}}\log(2n/m)} \\
    &= \frac{\Lambda(m)}{em^{iT}}\sum_{0 < r \leq \varepsilon m-1} \Big(\frac{(m+r)^{iT}}{r} + O\Big(\frac{\varepsilon}{r\log m}\Big)\Big) \\
    &= \frac{\Lambda(m)}{e}\sum_{1 \leq r \leq \varepsilon m-1} \frac{(1+\frac{r}{m})^{iT}}{r} +O(\varepsilon\log m). 
\end{align} 
We now estimate
\begin{equation} 
J(T) := \sum_{T^{1/4}\le r \le m\varepsilon}
\frac{1}{r}\exp\!\left(iT\log\left(1+\frac{r}{m}\right)\right) 
\end{equation}
where $\frac{T}{(1+\varepsilon)\log^{2}T} \leq m \leq T(1+\varepsilon)$. Let $\phi(r)=\frac{T}{2\pi}\log(1+\frac{r}{m})$ thus $\phi'(r)=\frac{T}{2\pi(r+m)}.$ Let $\eta \in (0, 1/2)$ be a constant, $\beta=\phi'(T^{1/4})$ and $\alpha=\phi'(m\varepsilon)$. By Lemma 4.10 of Titchmarsh [7, p.76], we have 
\begin{equation}
    J(T) = \sum_{\alpha-\eta \leq v \leq \beta+\eta} \int_{T^{1/4}}^{m\varepsilon} \frac{1}{x}e^{2\pi i (\phi(x)-vx)} \mathrm{d}x + O(T^{-1/4}\log T).
\end{equation} 
Since $\beta-\alpha \ll \varepsilon \log^{2} T = o(1)$, note that the range $[\alpha-\eta, \beta+\eta]$ contains at most one integer $v$. Summing (19) over $\frac{T}{(1+\varepsilon)\log^{2}T} \leq m \leq T(1+\varepsilon)$ gives
\begin{equation}
\sum_{\frac{T}{(1+\varepsilon)\log^{2}T} \leq m \leq T(1+\varepsilon)} J(T)
\ll \int_{T^{1/4}}^{T(1+\varepsilon)\varepsilon} \frac{1}{x} 
\Bigg| \sum_{m=\max\left(\frac{T}{(1+\varepsilon)\log^{2}T},\frac{x}{\varepsilon}\right)}^{T(1+\varepsilon)}
e^{2\pi i\big(\phi(x)-v x\big)} \Bigg| \, dx \;+\; T^{3/4}\log T,
\end{equation}
where $v=v(m)$ is the integer nearest to $\frac{T}{2\pi m}$. The dependence of $v$ on $m$ must be respected, so we write $k=v(m)$.  The condition $|k-\frac{T}{2\pi m}|\le\frac12$ implies
$m\in I_k:=\Bigl[\frac{T}{2\pi(k+\frac12)},\frac{T}{2\pi(k-\frac12)}\Bigr].$ 
For $m\ge T/((1+\varepsilon)\log^2 T)$ we have $0\leq k\ll\log^2 T$; hence $k$ runs over
$O(\log^2 T)$ nonnegative integers.  The length of $I_k$ is
$|I_k|\asymp T(k+1)^{-2}$. On $I_k$ the phase function is 
$g_k(m)=\frac{T}{2\pi}\log\!\Bigl(1+\frac{x}{m}\Bigr)-kx$ where $k$ is fixed. Its derivative satisfies $g_k'(m)=-\frac{T x}{2\pi m(m+x)}\asymp -\frac{T x}{m^{2}}$.
Applying the first‑derivative bound of van der Corput (Lemma 4.19 of Titchmarsh [7]) on each
interval $I_k$ yields
\[
\Bigl|\sum_{m\in I_k}e^{2\pi i g_k(m)}\Bigr|
\ll \frac{1}{\lambda_1}+|I_k|\lambda_1,
\qquad \lambda_1:=\frac{T x}{m^{2}}\asymp\frac{x(k+1)^2}{T}.
\]
Since $|I_k|\asymp \frac{T}{(1+k)^2}$, we have
\[
\Bigl|\sum_{m\in I_k}e^{2\pi i g_k(m)}\Bigr|
\ll \frac{T}{x}+x. 
\]
Summing over all admissible $k$ where $0\le k\ll\log^{2} T$, gives
\begin{equation}
\Bigl|\sum_{m} e^{2\pi i(\phi(x)-\nu x)}\Bigr|
\ll \sum_{k}\Bigl(\frac{T}{x}+x\Bigr)
\ll \frac{T\log^{2} T}{x}+x\log^{2}T.
\end{equation}
Inserting this bound into the integral in (20) gives 
\[
\int_{T^{1/4}}^{T(1+\varepsilon)\varepsilon} \frac{1}{x}
\Bigl|\sum_{m} e^{2\pi i(\phi(x)-\nu x)}\Bigr|\,dx
\ll \int_{T^{1/4}}^{T^{3/4}}\Bigl(\frac{T\log^{2} T}{x^{2}}+ \log^{2} T\Bigr)dx \ll T^{3/4}\log^{2} T .
\]
Plugging this into (20) gives 
\begin{equation}
\sum_{\frac{T}{(1+\varepsilon)\log^{2}T}\le m\le T(1+\varepsilon)}J(T)
\ll T^{3/4}\log^{2} T .
\end{equation}
Combining (20), (21) and (22) yields 
\begin{equation} 
\sum_{\frac{T}{(1+\varepsilon)\log^{2}T} \leq m \leq T(1+\varepsilon)}  \sum_{T^{1/4}\le r \le m\varepsilon}
\frac{1}{r}\exp\!\left(iT\log\left(1+\frac{r}{m}\right)\right) \ll T^{3/4}\log^{2} T.  
\end{equation}
  For 
  \begin{equation} 
  1\leq r \leq T^{1/4}
  \end{equation} we have $\log(1+\frac{r}{m})=\frac{r}{m} + O((r/m)^2)$ and $e^{O({Tr^{2}/m^2})}=1+O({Tr^{2}/m^2})$. Applying this together with (23) into the right-hand side of (17) yields 
\begin{align}
\sum_{\frac{T}{(1+\varepsilon)\log^{2}T} \leq m \leq T(1+\varepsilon)} \sum_{m+1 \leq 2n \leq m(1+\varepsilon)} \frac{\Lambda(m)(2n)^{iT}}{m^{1+\frac{1}{\log 2n}+iT}\log(2n/m)} &= \sum_{\frac{T}{(1+\varepsilon)\log^{2}T} \leq m \leq T(1+\varepsilon)}\frac{\Lambda(m)}{e}\sum_{1\leq r \leq T^{1/4}} \frac{1}{r}e^{\frac{irT}{m}} \\ 
&+ O( T^{3/4}\log^{2} T ).
\end{align}
For $\frac{m}{1+\varepsilon} \leq 2n \leq m-1,$ put $2n=m-r$ where $ 1\leq r \leq m(\varepsilon + O(\varepsilon^2))$ so that $\log(2n/m)=-\frac{r}{m}(1 + O(r/m))$ hence $\frac{1}{\log(2n/m)}=-\frac{m}{r}(1+O(r/m))=-\frac{m}{r}+O(1)$. Thus by arguing as above, one arrives at 
\begin{align}
\sum_{\frac{T}{(1+\varepsilon)\log^{2}T} \leq m \leq T(1+\varepsilon)} \sum_{\frac{m}{1+\varepsilon} \leq 2n \leq m-1} \frac{\Lambda(m)(2n)^{iT}}{m^{1+\frac{1}{\log 2n}+iT}\log(2n/m)} &= -\sum_{\frac{T}{(1+\varepsilon)\log^{2}T} \leq m \leq T(1+\varepsilon)}\frac{\Lambda(m)}{e}\sum_{1\leq r \leq T^{1/4}} \frac{1}{r}e^{-\frac{irT}{m}} \\ 
&+ O( T^{3/4}\log^{2} T).
\end{align}
Combining (28) and (26) gives 
\begin{align}
    &\sum_{\frac{T}{(1+\varepsilon)\log^{2}T}\leq m \leq T(1+\varepsilon)} \Big( \sum_{\frac{m}{1+\varepsilon} \leq 2n \leq m-1} + \sum_{m+1 \leq 2n \leq m(1+\varepsilon)}\Big) \frac{\Lambda(m)(2n)^{iT}}{m^{1+\frac{1}{\log 2n}+iT}\log(2n/m)} \\
    &= 
\frac{1}{e}\sum_{\frac{T}{(1+\varepsilon)\log^{2}T} \leq m \leq T(1+\varepsilon)}\Lambda(m)\sum_{1\leq r \leq T^{1/4}} \frac{e^{\frac{irT}{m}}-e^{-\frac{irT}{m}}}{r} + O(T^{3/4}\log^{2} T). 
\end{align}
Let $\lbrace{y\rbrace}$ denote the fractional part of $y$ and $\mid \mid y \mid \mid$ denote the distance to the integer nearest to $y$. Note that if $u>0$ is not an integer, we have 
\begin{equation}
    \sum_{n \leq X} \frac{\sin(2\pi nu)}{n}=\pi\Big(\frac{1}{2}-\lbrace{u\rbrace}  + O\Big(\frac{1}{1+X||u||} \Big)\Big)
\end{equation} as $X \rightarrow \infty$.  \footnote{see e.g. https://math.stackexchange.com/questions/3331510/error-estimation-for-the-fourier-series-of-the-fractional-part-of-x}
 Hence if $T \in \mathbb{Q}$ so that $T/2\pi m$ is irrational for $m \in \mathbb{N}$, we have
\begin{equation}
    \sum_{1\leq r \leq T^{1/4}} \frac{e^{\frac{irT}{m}}-e^{-\frac{irT}{m}}}{r} = \sum_{1\leq r \leq T^{1/4}} \frac{2i\sin(rT/m)}{r}= 2\pi i\Big(\frac{1}{2}-\lbrace{T/2\pi m \rbrace} + O \Big( \frac{1}{1+T^{1/4}||T/2\pi m||}\Big)\Big).
\end{equation} Define
\[
U(T)=\sum_{V\le m\le T(1+\varepsilon)}\frac{1}{1+T^{1/4}\,\Bigl\|\frac{T}{2\pi m}\Bigr\|},
\qquad V:=\frac{T}{(1+\varepsilon)\log^{2} T}.
\]
Set $Z=T^{1/4}$ and $x_m:=\frac{T}{2\pi m}.$ We split the sum $U(T)$ according to the size of $\|x_m\|$. For $\|x_m\|\le Z^{-1}$, we have $\frac{1}{1+Z\|x_m\|}\asymp 1,$
and the corresponding contribution is bounded by
\[
\#\{m\in[V,T(1+\varepsilon)]:\|x_m\|\le Z^{-1}\}.
\]
Fix $k\in\mathbb N$. The condition
$\Bigl|\frac{T}{2\pi m}-k\Bigr|\le Z^{-1}$ 
is equivalent to $\frac{T}{2\pi(k+Z^{-1})}
\le m\le
\frac{T}{2\pi(k-Z^{-1})}.$ 
The length of this interval is $\frac{T}{2\pi}
\left(\frac{1}{k-Z^{-1}}-\frac{1}{k+Z^{-1}}\right)
=
\frac{T}{2\pi}\cdot
\frac{2Z^{-1}}{k^2-Z^{-2}}
\ll
\frac{T}{Zk^2}$. 
Hence, for fixed $k$, the number of such integers $m$ is $\ll \frac{T}{Z k^2}+1.$
Since $x_m=\frac{T}{2\pi m}\in
\Bigl[\frac{1}{2\pi (1+\varepsilon)},\frac{(1+\varepsilon)\log^{2} T}{2\pi}\Bigr],$ 
we have $1\le k\ll(\log T)^2$. Summing over $k$ gives
\[\#\{m:\|x_m\|\le Z^{-1} \} \ll \sum_{k\le(\log T)^2} \left(\frac{T}{Z k^2}+1\right)
\ll \frac{T}{Z}.\]
Thus the contribution from the $x_m$ for which $\|x_m\|\le Z^{-1}$, is 
\begin{equation} 
\sum_{\|x_m\|\le Z^{-1}}
\frac{1}{1+Z\|x_m\|} \ll \frac{T}{Z} =T^{3/4}.
\end{equation}
For $\|x_m\|>Z^{-1}$ we use $\frac{1}{1+Z\|x_m\|}\le\frac{1}{Z\|x_m\|}.$ Partition $(0,1/2]$ dyadically: $2^{-j-1}<\|x_m\|\le2^{-j}$ where $0\le j\le \log_2 Z.$
For such $m$, $\frac{1}{1+Z\|x_m\|}\ll\frac{2^j}{Z}$. Since $x_m$ is monotone in $m$, the number of $m$ in each dyadic interval is $\ll T\cdot2^{-j}.$
Hence the contribution of each $j$ is $\ll T\cdot2^{-j}\cdot\frac{2^j}{Z}
=
\frac{T}{Z}.$ 
Summing over $0\leq j\le\log_{2}Z$ yields 
\begin{equation}
\sum_{\|x_m\|>Z^{-1}} \frac{1}{1+Z\|x_m\|} \ll \frac{T}{Z}\log Z \ll T^{3/4}\log T.
\end{equation}
Combining both ranges gives
\begin{equation}
U(T)=\sum_{\frac{T}{(1+\varepsilon)\log^{2} T}\le m\le T(1+\varepsilon)}
\frac{1}{1+T^{1/4}\Bigl\|\tfrac{T}{2\pi m}\Bigr\|} \ll T^{3/4}\log T.
\end{equation}
Inserting (32) into (30) and applying (35) gives
\begin{align}
    &\sum_{\frac{T}{(1+\varepsilon)\log^{2}T}\leq m \leq T(1+\varepsilon)} \Big( \sum_{\frac{m}{1+\varepsilon} \leq 2n \leq m-1} + \sum_{m+1 \leq 2n \leq m(1+\varepsilon)}\Big) \frac{\Lambda(m)(2n)^{iT}}{m^{1+\frac{1}{\log 2n}+iT}\log(2n/m)} \\
&= \frac{2\pi i}{e}\sum_{\frac{T}{(1+\varepsilon)\log^{2}T} \leq m \leq T(1+\varepsilon)}\Lambda(m)\Big(\frac{1}{2}-\lbrace{T/2\pi m \rbrace}\Big) + O(T^{3/4}\log^{2} T). 
\end{align} 
Note that $\sum_{m \leq x} \Lambda(m)\lbrace{x/m \rbrace} = (1-B)x + o(x)$ where $B=0.57721
\cdots$ denotes the Euler-Mascheroni constant [5, p.249]. Also note that $\lbrace{\frac{T}{2\pi m} \rbrace}=\frac{T}{2\pi m} $ for $\frac{T}{2\pi} < m \leq T(1+\varepsilon)$. Thus 
\begin{align}
    & \sum_{m \leq \frac{T}{(1+\varepsilon)\log^{2} T}} \Lambda(m)\lbrace{T/2\pi m\rbrace} - \sum_{m \leq \frac{T}{2\pi }} \Lambda(m)\lbrace{T/2\pi m \rbrace} - \sum_{\frac{T}{2\pi} < m \leq T(1+\varepsilon)} \Lambda(m)\lbrace{T/2\pi m \rbrace} \\
    &+ \sum_{\frac{T}{(1+\varepsilon)\log^{2}T} \leq m \leq T(1+\varepsilon)} \frac{\Lambda(m)}{2} \\
    &=(B-1)\frac{T}{2\pi}  - \frac{T}{2\pi}\sum_{\frac{T}{2\pi}<m\leq T(1+\varepsilon)} \frac{\Lambda(m)}{m} + \frac{T}{2} + o(T) \\
    &= (B-1)\frac{T}{2\pi} - \frac{T}{2\pi}(\log T -\log (T/2\pi)+O(\varepsilon)) + \frac{T}{2} + o(T) \\
    &= (B-1)\frac{T}{2\pi} - \frac{T}{2\pi}\log (2\pi)+\frac{T}{2} + o(T)
\end{align} 
since $\sum_{m \leq x} \Lambda(m) \sim x$ and $\sum_{m \leq x} \frac{\Lambda(m)}{m} = \log x -B +o(1)$ [5, p.187]. Plugging (42) into (37) gives 
\begin{align}
    &\sum_{\frac{T}{(1+\varepsilon)\log^{2}T}\leq m \leq T(1+\varepsilon)} \Big( \sum_{\frac{m}{1+\varepsilon} \leq 2n \leq m-1} + \sum_{m+1 \leq 2n \leq m(1+\varepsilon)}\Big) \frac{\Lambda(m)(2n)^{iT}}{m^{1+\frac{1}{\log 2n}+iT}\log(2n/m)} \\
&= \frac{2\pi i}{e}\Big((B-1)\frac{T}{2\pi} - \frac{T}{2\pi}\log (2\pi) + \frac{T}{2} + o(T) \Big). 
\end{align} 
By a similar argument, one arrives at 
\begin{align}
    &\sum_{\frac{T}{(1+\varepsilon)\log^{2}T}\leq m \leq T(1+\varepsilon)} \Big( \sum_{\frac{m}{1+\varepsilon} \leq 2n \leq m-1} + \sum_{m+1 \leq 2n \leq m(1+\varepsilon)}\Big) \frac{\Lambda(m)(2n)^{iH}}{m^{1+\frac{1}{\log 2n}+iH}\log(2n/m)} \\
&= \frac{2\pi i}{e}\Big((B-1)\frac{H}{2\pi} - \frac{H}{2\pi}\log \pi +\frac{T}{2} + o(T) \Big) \\
&= \frac{2\pi i}{e}\Big((B-1)\frac{T}{\pi} - \frac{T}{\pi}\log \pi +\frac{T}{2} + o(T) \Big)
\end{align} 
since $H=2T+O(1)$. Combining (47) and (44) gives 
\begin{align}
  & \sum_{\frac{T}{(1+\varepsilon)\log^{2}T}\leq m \leq T(1+\varepsilon)} \Big( \sum_{\frac{m}{1+\varepsilon} \leq 2n \leq m-1} + \sum_{m+1 \leq 2n \leq m(1+\varepsilon)}\Big) \frac{\Lambda(m)(2n)^{iT}}{m^{1+\frac{1}{\log 2n}+iT}\log(2n/m)}  \\
  & - \sum_{\frac{T}{(1+\varepsilon)\log^{2}T}\leq m \leq T(1+\varepsilon)} \Big( \sum_{\frac{m}{1+\varepsilon} \leq 2n \leq m-1} + \sum_{m+1 \leq 2n \leq m(1+\varepsilon)}\Big) \frac{\Lambda(m)(2n)^{iH}}{m^{1+\frac{1}{\log 2n}+iH}\log(2n/m)} \\
  &= \frac{2\pi i}{e}\Big((1-B)\frac{T}{2\pi} + \frac{T}{2\pi}\log (\pi/2) + o(T) \Big).
\end{align}
Note that $\sum_{m \leq x}\frac{\Lambda(2m)}{2m} = \Big(\sum_{m=1}^{\infty} -\sum_{m \geq x+1} \Big) \frac{\Lambda(2m)}{2m} =\log 2 + O(1/x)$. Therefore, summing (9) for $\frac{T}{\log^{2}T}<2n\leq T$, changing the order  of summation, invoking (50), (14), (13) and (12) yields
\begin{align}
    \sum_{\frac{T}{\log^{2} T} < 2n \leq T} I_1 &= \frac{T-H}{2\pi }\sum_{\frac{T}{\log^{2}T}<2n\leq T} \frac{\Lambda(2n)}{2n} + (1-B)\frac{T}{2\pi} + \frac{T}{2\pi}\log (\pi/2) + o(T) \\
    &=  (1-B)\frac{T}{2\pi} + \frac{T}{2\pi}\log (\pi/2) + o(T).
\end{align}
By Lemma 2, we have  
\begin{equation}
\sum_{T \leq \gamma_{T} \leq H} (2n)^{\rho_{T}-1}  =-\frac{T\Lambda(2n)}{4\pi n}  + O((\log T)\log \log 6n) 
\end{equation} 
uniformly for $1 \leq n \ll T$. Hence, summing (53) for $1\leq n \leq \frac{T}{2\log^{2} T}$ gives 
\begin{align}
    \sum_{1 \leq n \leq \frac{T}{2\log^{2} T}} \sum_{T \leq \gamma_{T} \leq H} (2n)^{\rho_{T}-1}   &= -\frac{T}{2\pi}\sum_{1 \leq n \leq \frac{T}{2\log^{2} T}} \frac{\Lambda(2n)}{2n} + o(T) \\
    &=-\frac{T}{2\pi}\log 2 + o(T). 
\end{align} 
Combining (55), (52), (8) and (3) yields 
\begin{align}
   \sum_{1 \leq n \leq \frac{T}{2}} \sum_{T \leq \gamma_{T} \leq H} (2n)^{\rho_{T}-1} &= \Big(\sum_{1\leq n \leq \frac{T}{2\log^{2} T}} + \sum_{\frac{T}{2\log^{2}T} < n \leq \frac{T}{2}} \Big) \sum_{T \leq \gamma_{T} \leq H} (2n)^{\rho_{T}-1} \\
    &= (1-B+\log (\pi/4))\frac{T}{2\pi}  + o(T)
\end{align} 
\textbf{for the sequence of rational $T, H \in [2T, 2T+1]$ that satisfy (1).} Suppose that $\Theta<1$. Then by Lemma 1 we have 
\begin{equation}
    \sum_{1\leq n \leq \frac{T}{2}}(2n)^{\rho_{T}-1} =2^{\rho_{T}-1}\sum_{1\leq n \leq \frac{T}{2}} n^{\rho_{T}-1} \ll T^{\Theta-1}
\end{equation}
\textbf{uniformly for} $T \leq \gamma_{T} \leq H.$  For large $x>0$, note that there are $\asymp x \log x$ zeros of $\zeta$ with imaginary parts in $(0, x]$ [5, Corollary 14.3]. Hence 
\begin{equation}
    \sum_{1\leq n \leq \frac{T}{2}} \sum_{T \leq \gamma_{T} \leq H} (2n)^{\rho_{T}-1} = \sum_{T \leq \gamma_{T} \leq H}  \sum_{1\leq n \leq \frac{T}{2}} (2n)^{\rho_{T}-1}  \ll T^{\Theta}\log T
\end{equation} 
\textbf{uniformly} for $T \geq T_0$. But (59) is a contradiction to (57), thus our supposition must be false. This completes the proof. 
\end{proof} 
\textbf{Corollary 1.} \textit{There exists an infinite sequence of zeros of $\zeta$, whose real parts converge to 1.} 
\begin{proof}  
Combine Theorem 1 with the fact that $\zeta(s) \neq 0$ for $\Re(s) \geq 1$ [5, Theorem 6.6].
\end{proof}
\textbf{Remark 1.} Notice that our argument actually shows that, \textit{for all large enough $T \in \mathbb{Q}^+$, there exists some $\rho_T$ with $T \leq \gamma_{T} \leq 2T+1$ and $\Re(\rho_T) = 1 - o(1)$.}
\\
\\
\textbf{Remark 2.} Our argument crucially relies on: 
(i) the functional equation for $\zeta$, and  
(ii) theorem 4.11 of [7]. Both (i) and (ii) do not hold for Beurling zeta functions corresponding to Beurling number systems whose integers are not the ordinary integers. Hence our approach wouldn't work for such, whose analogues of the RH are known to be true. Also, (ii) does not hold for the Weil zeta functions, which are known to satisfy the RH. 
\\
\\
\textbf{Remark 3.} From the literature on arXiv and elsewhere, it doesn't appear like this approach on RH has been tried before. However, the idea of applying Gonek's formula (53) in estimating sums over zeros of zeta in general, is not new. As one can confirm from the literature, it has proven to be a useful tool for such tasks. It is therefore instructive to sum (9) carefully, e.g. by not rushing to naively bound $(2n)^{it}$ by 1 for $2n$ not too close to $m$. Doing so will definitely make one lose important information needed for a sufficiently sharp bound for the sums over $2n \leq \frac{m}{1+\varepsilon}$ and $2n \geq m(1+\varepsilon)$. Note that the analogue of (53) from the Riemann explicit formula is $\sum_{T \leq \gamma_{T} \leq H} \frac{(2n)^{\rho_{T}-1}}{\rho_T} \ll \frac{(\log T)(\log \log T)}{T}.$ This has no main term, so it's not quite as elegant/useful as (53). 
\\
\\
\textbf{Remark 4.} The reader might wonder if our approach can give a lower bound for $\Theta_T:$ the supremum of the real parts of the zeros $\rho_T$ with $T \leq \gamma_T \leq H.$ For that, one would need our Lemma 1 to also hold for $\sigma=\sigma_t = o(1)$ with an error term $\ll (\log t)^{-1}$. However, this is not possible (as can be seen e.g by estimating $\zeta(a+it)=\sum_{n \leq \frac{T}{2}} n^{-a-it} + O(T^{-a})=\sum_{n \leq T/2} n^{(-a+\sigma_{t})-\sigma_t-it}+O(T^{-a})$ for any fixed $a>0$ via Abel summation, under the assumption that $\sum_{n \leq T/2} n^{-\sigma_{t}-it}=\zeta(\sigma_t + it) + O(1/\log t)$). It can also be seen that Titchmarsh's proof of his Theorem 4.11 doesn't hold for $\sigma = \sigma_t = o(1)$. Indeed, he begins his argument by quoting the Euler-McLaurin summation formula for $\zeta(\sigma+it)$ (equation 3.5.3). This is initially defined for $\sigma >1$, and then is extended by holomorphic continuation to $\sigma>0$. Recall that a function can only be holomorphically extended on open sets, and that's precisely why Titchmarsh assumes that $\sigma$ is bounded below by some constant $\sigma_0 > 0$. For the same reason, one \textbf{cannot} apply any of the best known zero-free regions of  $\zeta$ (e.g. Vinogradov-Koborov) to bound $\Re(\rho_{T})$ in (58).
\\

\end{document}